\documentclass{amsart}
\usepackage{latexsym}
\usepackage{amssymb}
\usepackage{amsmath}

\begin{document}


\newtheorem{thm}{Theorem}
\newtheorem{problem}{Problem}
\newtheorem{definition}{Definition}
\newtheorem{lemma}{Lemma}
\newtheorem{proposition}{Proposition}
\newtheorem{corollary}{Corollary}
\newtheorem{example}{Example}
\newtheorem{conjecture}{Conjecture}
\newtheorem{algorithm}{Algorithm}
\newtheorem{exercise}{Exercise}
\newtheorem{xample}{Example}
\newtheorem{remarkk}{Remark}

\newcommand{\be}{\begin{equation}}
\newcommand{\ee}{\end{equation}}
\newcommand{\bea}{\begin{eqnarray}}
\newcommand{\eea}{\end{eqnarray}}
\newcommand{\beq}[1]{\begin{equation}\label{#1}}
\newcommand{\eeq}{\end{equation}}
\newcommand{\beqn}[1]{\begin{eqnarray}\label{#1}}
\newcommand{\eeqn}{\end{eqnarray}}
\newcommand{\beaa}{\begin{eqnarray*}}
\newcommand{\eeaa}{\end{eqnarray*}}
\newcommand{\req}[1]{(\ref{#1})}

\newcommand{\lip}{\langle}
\newcommand{\rip}{\rangle}
\newcommand{\uu}{\underline}
\newcommand{\oo}{\overline}
\newcommand{\La}{\Lambda}
\newcommand{\la}{\lambda}
\newcommand{\eps}{\varepsilon}
\newcommand{\om}{\omega}
\newcommand{\Om}{\Omega}
\newcommand{\ga}{\gamma}
\newcommand{\ka}{\kappa}
\newcommand{\rrr}{{\Bigr)}}
\newcommand{\qqq}{{\Bigl\|}}

\newcommand{\dint}{\displaystyle\int}
\newcommand{\dsum}{\displaystyle\sum}
\newcommand{\dfr}{\displaystyle\frac}
\newcommand{\bige}{\mbox{\Large\it e}}
\newcommand{\integers}{{\Bbb Z}}
\newcommand{\rationals}{{\Bbb Q}}
\newcommand{\reals}{{\rm I\!R}}
\newcommand{\realsd}{\reals^d}
\newcommand{\realsn}{\reals^n}
\newcommand{\NN}{{\rm I\!N}}
\newcommand{\DD}{{\rm I\!D}}
\newcommand{\degree}{{\scriptscriptstyle \circ }}
\newcommand{\dfn}{\stackrel{\triangle}{=}}
\def\complex{\mathop{\raise .45ex\hbox{${\bf\scriptstyle{|}}$}
     \kern -0.40em {\rm \textstyle{C}}}\nolimits}
\def\hilbert{\mathop{\raise .21ex\hbox{$\bigcirc$}}\kern -1.005em {\rm\textstyle{H}}} 
\newcommand{\RAISE}{{\:\raisebox{.6ex}{$\scriptstyle{>}$}\raisebox{-.3ex}
           {$\scriptstyle{\!\!\!\!\!<}\:$}}} 

\newcommand{\hh}{{\:\raisebox{1.8ex}{$\scriptstyle{\degree}$}\raisebox{.0ex}
           {$\textstyle{\!\!\!\! H}$}}}

\newcommand{\OO}{\won}
\newcommand{\calA}{{\mathcal A}}
\newcommand{\calB}{{\mathcal B}}
\newcommand{\calC}{{\cal C}}
\newcommand{\calD}{{\cal D}}
\newcommand{\calE}{{\mathcal E}}
\newcommand{\calF}{{\mathcal F}}
\newcommand{\calG}{{\cal G}}
\newcommand{\calH}{{\cal H}}
\newcommand{\calK}{{\cal K}}
\newcommand{\calL}{{\mathcal L}}
\newcommand{\calM}{{\mathcal M}}
\newcommand{\calO}{{\cal O}}
\newcommand{\calP}{{\cal P}}
\newcommand{\calT}{{\mathcal T}} 
\newcommand{\calU}{{\mathcal U}}
\newcommand{\calX}{{\cal X}}
\newcommand{\calY}{{\mathcal Y}}
\newcommand{\calZ}{{\mathcal Z}}
\newcommand{\calXX}{{\cal X\mbox{\raisebox{.3ex}{$\!\!\!\!\!-$}}}}
\newcommand{\calXXX}{{\cal X\!\!\!\!\!-}}
\newcommand{\gi}{{\raisebox{.0ex}{$\scriptscriptstyle{\cal X}$}
\raisebox{.1ex} {$\scriptstyle{\!\!\!\!-}\:$}}}
\newcommand{\intsim}{\int_0^1\!\!\!\!\!\!\!\!\!\sim}
\newcommand{\intsimt}{\int_0^t\!\!\!\!\!\!\!\!\!\sim}
\newcommand{\pp}{{\partial}}
\newcommand{\al}{{\alpha}}
\newcommand{\sB}{{\cal B}}
\newcommand{\sL}{{\cal L}}
\newcommand{\sF}{{\cal F}}
\newcommand{\sE}{{\cal E}}
\newcommand{\sX}{{\cal X}}
\newcommand{\R}{{\rm I\!R}}
\renewcommand{\L}{{\rm I\!L}}
\newcommand{\vp}{\varphi}
\newcommand{\N}{{\rm I\!N}}
\def\ooo{\lip}
\def\ccc{\rip}
\newcommand{\ot}{\hat\otimes}
\newcommand{\rP}{{\Bbb P}}
\newcommand{\bfcdot}{{\mbox{\boldmath$\cdot$}}}

\renewcommand{\varrho}{{\ell}}
\newcommand{\dett}{{\textstyle{\det_2}}}
\newcommand{\sign}{{\mbox{\rm sign}}}
\newcommand{\TE}{{\rm TE}}
\newcommand{\TA}{{\rm TA}}
\newcommand{\E}{{\rm E\,}}
\newcommand{\won}{{\mbox{\bf 1}}}
\newcommand{\Lebn}{{\rm Leb}_n}
\newcommand{\Prob}{{\rm Prob\,}}
\newcommand{\sinc}{{\rm sinc\,}}
\newcommand{\ctg}{{\rm ctg\,}}
\newcommand{\loc}{{\rm loc}}
\newcommand{\trace}{{\,\,\rm trace\,\,}}
\newcommand{\Dom}{{\rm Dom}}
\newcommand{\ifff}{\mbox{\ if and only if\ }}
\newcommand{\nproof}{\noindent {\bf Proof:\ }}
\newcommand{\remark}{\noindent {\bf Remark:\ }}
\newcommand{\remarks}{\noindent {\bf Remarks:\ }}
\newcommand{\note}{\noindent {\bf Note:\ }}

\newcommand{\boldx}{{\bf x}}
\newcommand{\boldX}{{\bf X}}
\newcommand{\boldy}{{\bf y}}
\newcommand{\boldR}{{\bf R}}
\newcommand{\uux}{\uu{x}}
\newcommand{\uuY}{\uu{Y}}

\newcommand{\limn}{\lim_{n \rightarrow \infty}}
\newcommand{\limN}{\lim_{N \rightarrow \infty}}
\newcommand{\limr}{\lim_{r \rightarrow \infty}}
\newcommand{\limd}{\lim_{\delta \rightarrow \infty}}
\newcommand{\limM}{\lim_{M \rightarrow \infty}}
\newcommand{\limsupn}{\limsup_{n \rightarrow \infty}}

\newcommand{\ra}{ \rightarrow }

\newcommand{\ARROW}[1]
  {\begin{array}[t]{c}  \longrightarrow \\[-0.2cm] \textstyle{#1} \end{array} }

\newcommand{\AR}
 {\begin{array}[t]{c}
  \longrightarrow \\[-0.3cm]
  \scriptstyle {n\rightarrow \infty}
  \end{array}}

\newcommand{\pile}[2]
  {\left( \begin{array}{c}  {#1}\\[-0.2cm] {#2} \end{array} \right) }

\newcommand{\floor}[1]{\left\lfloor #1 \right\rfloor}

\newcommand{\mmbox}[1]{\mbox{\scriptsize{#1}}}

\newcommand{\ffrac}[2]
  {\left( \frac{#1}{#2} \right)}

\newcommand{\one}{\frac{1}{n}\:}
\newcommand{\half}{\frac{1}{2}\:}

\def\le{\leq}
\def\ge{\geq}
\def\lt{<}
\def\gt{>}

\def\squarebox#1{\hbox to #1{\hfill\vbox to #1{\vfill}}}
\newcommand{\nqed}{\hspace*{\fill}
          \vbox{\hrule\hbox{\vrule\squarebox{.667em}\vrule}\hrule}\bigskip}
\title{Log-concave Measures }
\author{D. Feyel and A. S. \"Ust\"unel}
\address{D.F.: 10, Quai des Platreries, 77920 Samois/Seine, France}
\email{denis.feyel@orange.fr}
\address{A.S.U.: Institut Telecom-Paristech, 46, rue Barrault, 75013 Paris, France }
\email{ustunel@enst.fr}

\subjclass{Primary 60B, 60H; Secondary 28C}
\keywords{log-concave measures, logarithmic Sobolev inequality,
  measure transport, Pr\'kopa-Leindler inequality, $H$-convexity}

\maketitle
\begin{abstract}
\noindent
We study  the log-concave
measures, their characterization via the Pr\'ekopa-Leindler property
and also define a subset of it whose elements are called
super log-concave measures which have the property of satisfying a
logarithmic Sobolev inequality. We give some results about their
stability. Certain relations with measure transportation are also indicated.
\end{abstract}

\section{Introduction}
\noindent
The importance of logarithmically concave functions and measures has
been discovered in the 70's
(cf. \cite{Prekopa,Prekopa-1,Leindler,Bor-1}) and they have found
applications  immediately in physics (cf. \cite{Simon}).  This notion has
gained further importance when its  close relations to Monge-Amp\`ere
equation and more generally to measure transportation has been
realized (cf.\cite{Vil} and the references there). In this work we
give a general treatment of the subject beginning from the finite
dimensional case and going towards to the infinite dimensions. A
preliminary definition of log-concave measures is given using the
Lebesgue measure then we extend this definition to the measures
without density by using the Pr\'ekopa-Leindler  property which
characterizes them, which is also equivalent to Brunn-Minkowski
property. We also introduce the notion of super log-concave measures,
namely these are the measures which decrease very rapidly at infinity,
in fact faster than some Gaussian measure and they always satisfy the
logarithmic Sobolev inequality. Their definition uses the Euclidean
structure of the underlying space on which they are defined,
consequently, we need a special structure if we want to extend this
notion to the infinite dimensional case. This is done by adjoining a
rigged Hilbert space structure  to the Fr\'echet space supporting the
measure under question, whose typical example is an abstract Wiener
space with its Cameron-Martin space. Afterwards we concentrate
ourselves to the case of Wiener space, first we give some
complementary results about the Jacobians corresponding to the image
of the Wiener measure under a  general perturbation of identity with
lacking regularity and show that although each term of the Jacobian
is not properly defined, their multiplication may create a
renormalization, which is the typical case with monotone shifts. Even
in this case, the log-concave character of the Jacobian is preserved
and we use this observation to give another proof of the
Pr\'ekopa-Leindler property for the Gaussian measure with ``less
log-concave'' functions (called $1$-log concave) and we give also another proof of this
property using the reverse martingale convergence theorem, in
particular we prove that the property of the functions used to test
Pr\'ekopa-Leindler property is preserved under the conditional
expectations and the Ornstein-Uhlenbeck semigroup.

\section{Preliminaries and notations}
Generally each section contains its notational conventions with the
exception of the Gaussian case; in fact we denote by $(W,H,\mu)$ an
abstract Wiener space, namely,  $W$ is a separable Fr\'echet space,
$H\subset W$ is a separable Hilbert space densely injected into $W$,  $\mu$
is the unit Gauss measure supported by $W$ which is quasi-invariant
under the translation by the elements of $H$. One defines the usual
Gateaux  derivative of the nice functions on $W$ along the subspace
$H$, due to the quasi-invariance, this derivative has a unique
(closed) extension to all the $L^p(\mu)$-spaces with $p\geq 1$ and it
is denoted by $\nabla$ and called Sobolev derivative
(cf. \cite{ASU} for instance). Hence, for
a nice function $f$  on $W$, $\nabla f$ defines a linear functional on
$H$ $\mu$-almost surely, consequently it can be identified with an
element of $H^\star=H$. Since the Sobolev derivative maps the scalar
functions to $H$-valued functions, its adjoint,  called the divergence
operator, maps the vector valued
functions to scalar ones and denoted by $\delta$, in particular, we
have 
$$
\int_W\,(\xi,\nabla f)_Hd\mu=\int_Wf\,\delta \xi d\mu\,,
$$
for $\xi:W\to H$ cylindrical, well-known as the integration by parts
formula. The very remarkable property of the divergence is that in the
case of classical Wiener space $W=C([0,1],\R^n)$, if $\xi$ has an
adapted Lebesgue density $\xi'(t,w)$, then it  holds that
$$
\delta\xi(w)=\int_0^1\xi'(t,w)dW_t(w)\,,
$$
where the integral with $dW$ denotes the It\^o stochastic integral.

\section{Log-concave and super log-concave measures  in finite dimension}
\noindent
We begin by an initial concept:
\begin{definition}
A (positive) measure $\rho$ on $\R^d$ is said to satisfy the
Pr\'ekopa-Leindler property if for any positive, continuous functions
of compact support, say $a,b,c$ such that 
\be
\label{PK-cond}
a(sx+ty)\geq b(x)^s c(y)^t
\ee
for any $x,y\in\R^d,\,s+t=1$, one has
\be
\label{P-L}
\rho(a)\geq \rho(b)^s\rho(c)^t\,.
\ee
\end{definition}
\noindent
\begin{thm}
\label{PK-thm}
Assume that $\theta$ is a  log-concave function, denote by $\rho$ the
measure $d\rho(x)=\theta(x)dx$. Then $\rho$ satisfies the Pr\'ekopa-Leindler property.
\end{thm}

\noindent
Here are  two other well-known results which are due to Pr\'ekopa,
\cite{Prekopa},  that we derive using the above considerations:
\begin{corollary}
\label{cor-1}
Let $f$ and $g$ be two integrable, log-concave functions, then their
convolution $f\star g$ is  log-concave. Moreover, for any
$c\in\R$, the function  $f\star_cg$ which is defined as
$$
f\star_cg(x)=\int f(cx-y)g(y)dy
$$
is again a log-concave function.
\end{corollary}
\nproof
Let $f_x(u)=f(x-u)$, then $(f_{sx+ty}g)(su+tv)\geq
  (f_x(u)g(u))^s(f_y(v)g(v))^t$, hence we can apply the theorem and
  obtain the inequality which characterizes the log-concavity. The
  proof of the second part is similar.
\nqed

\begin{corollary}
\label{cor-2}
Let $f(x,\xi)$ be a log-concave function on $\R^d\times \R^m$, define
$F(x)$ as to be
$$
F(x)=\int_{\R^m}f(x,\xi)d\xi\,.
$$
Then $F$ is log-concave on $\R^d$.
\end{corollary}
\nproof
Let $f_x(\xi)$ denote $f(x,\xi)$, then for
$s+t=1,\,x,y\in\R^d,\,\xi,\eta\in\R^m$, we have 
$$
f_{sx+ty}(s\xi+t\eta)\geq f_x(\xi)^sf_y(\eta)^t\,,
$$
hence 
\beaa
F(sx+ty)&=&\int f_{sx+ty}(\xi)d\xi\geq \left(\int f_x(\xi)d\xi\right)^s\left(\int
f_y(\xi)d\xi\right)^t\\
&=&F(x)^sF(y)^t\,.
\eeaa
\nqed

\noindent
The results above indicate that the  following definition is reasonable
\begin{definition}
\label{defn-1}
A (positive) measure on $\R^d$ is called log-concave if any of its convolutions
with log-concave continuous functions of compact support has a
log-concave density.
\end{definition}
\remark
Let us  note that, for a measure $\rho$ to be log-concave, it suffices the
existence of just one continuous, log-concave function $\theta$ such
that $\theta_\eps\star\rho$ has a log-concave density for any $\eps>0$
where $\theta_\eps$ denotes the $d$-dimensional rescaling of $\theta$. In fact 
using the commutativity of convolutions, we obtain also that
$\tilde{\theta}\star\rho$ has again a log-concave density, for any
other continuous log-concave function of compact support
$\tilde{\theta}$. 
\begin{proposition}
Let $\rho_1$ and $\rho_2$ be two log-concave measures on $\R^n$ and
$\R^m$ respectively. Then the product measure $\rho_1\otimes\rho_2$ is
log-concave on $\R^n\times \R^m$.
\end{proposition}
\nproof 
From Definition \ref{defn-1}, it suffices to assume that
$d\rho_1(x)=e^{-V_1(x)}dx$ and $d\rho_2(y)=e^{-V_2(y)}dy$, where
$V_i,\,i=1,2$ are convex functions. Then $(x,y)\to V_1(x)+V_2(y)$ is
convex on $\R^n\times \R^m$.
\nqed

\begin{proposition}
Let $(\rho_n,n\geq 1)$ be a sequence of log-concave measures
converging weakly to $\rho$, then $\rho$ is also log-concave.
\end{proposition}
\nproof 
Let $\theta$ be a log-concave function as described above, then the
density of $\theta\star\rho$ is the limit of the densities of
$(\theta\star\rho_n)$.
\nqed

\noindent
The following theorem is the first pavement to extend the definition
of log-concavity to the infinite dimensional case where the Lebesgue
measure does not exist:
\begin{thm}
\label{carac-thm}
  $\rho$ is a log-concave measure if and only if it satisfies the
  Pr\'ekopa-Leindler property.
\end{thm}
\nproof
Suppose that $\rho$ has a continuous density $F$ with respect to the
Lebesgue measure of $\R^d$. Let $f=1_A,\,g=1_B$, where $A$ and $B$ are
two balls whose centers are  located at $a$ and $b$ respectively and
whose diameter will tend to zero. Let $C=\half (A+B)$  and $h=1_C$ then
$h(\alpha x+\beta y)\geq f(x)^\alpha g(y)^\beta$ for any
$\alpha+\beta=1, \,x,y\in \R^d$. Hence for $\alpha=\half$, by taking
the limit, we obtain
$$
F(\half(a+b))\geq F(a)^\half F(b)^\half\,,
$$
which is a sufficient condition for the log-concavity of $F$. The
general case follows by taking the convolution of $\rho$. The
necessity follows from the transport argument that we have used in the
proof of Theorem \ref{PK-thm}.
\nqed

\begin{corollary}
If a measure satisfies the  relation (\ref{P-L})  for the continuous functions
satisfying the condition (\ref{PK-cond}), then it  also satisfies the same
relation  for Borel functions satisfying  (\ref{PK-cond}).
\end{corollary}
\nproof
If the relation (\ref{P-L}) is satisfied by $\rho$, then it is also
satisfied by $\rho\star\theta$ where $\theta$ is a log-concave
function,  since 
$\rho(\theta\star a)=(\rho\star\theta)(a)$.
\nqed
\begin{corollary}
\label{density-cor}
Let $\rho$ be a log-concave measure and let $F$ be a convex function,
then the measure $\nu$ defined as
$$
d\nu(x)=e^{-F(x)}d\rho(x)
$$
is again log-concave.
\end{corollary}
\nproof
The new measure obviously  satisfies the Pr\'ekopa-Leindler property.
\nqed

\noindent
A very close characterization of the log-concave measures can be given
by the Brunn-Minkowski inequality whose proof is similar to that of
Theorem \ref{carac-thm}
\begin{thm}
The measure $\rho$ is log-concave if and only if 
$$
\rho(sA+tB)\geq \rho(A)^s\rho(B)^t
$$
for any measurable $A,B$ and $s+t=1$.
\end{thm}

\begin{definition}
A (positive) measure $\rho$ on $\R^n$ is called $\alpha$-super log-concave
($\alpha$-s.l.c in short) if the measure
$$
e^{\frac{\alpha}{2}|x|^2}\rho(dx)
$$
is a log-concave measure, where $\alpha\geq 0$ and $|\cdot|$ denotes
the Euclidean norm.
\end{definition}
\begin{remarkk}
 Since $\exp-\frac{\alpha}{2}|x|^2$ is a log-concave function,
any $\alpha$-s.l.c. measure is log-concave.
\end{remarkk}
\begin{proposition}
Assume that a measure $\rho$ can be represented as
$$
d\rho(x)=e^{-V(x)}dx
$$
where $V$ is a locally integrable, lower bounded function such that
$$
\nabla^2V\geq \al I_{\R^n}
$$
in the sense of distributions, then $\rho$ is $\alpha$-s.l.c.
`\end{proposition}
\nproof
Evidently the condition implies the convexity of the function
$x\to V(x)-\frac{\al}{2}|x|^2$.
\nqed
\begin{remarkk}
If $\rho_i$ are finitely many  $\alpha_i$-s.l.c. measures on
$\R^{n_i}$, then their product is an $\min_i\alpha_i$-s.l.c. measure.
\end{remarkk}
\noindent
The proof of the following lemma follows from  that of Theorem
\ref{carac-thm}:
\begin{lemma}
A measure $\rho$ is an $\alpha$-s.l.c. if and only if for any $a,b,c$
continuous, positive functions of compact support such that, for any
$x,y\in \R^n$, $s+t=1$, 
$a(sx+ty)\geq b(x)^sc(y)^t$, one has 
$$
\rho(a_\alpha)\geq \rho(b_\alpha)^s\rho(c_\alpha)^t\,,
$$
where, for a given function $f$, $f_\alpha$ is defined as
$$
f_\alpha(x)=\exp(\frac{\alpha}{2}|x|^2)f(x)\,.
$$
\end{lemma}

\noindent
The proof of the following is obvious:
\begin{lemma}
\label{cvg-lemma}
Suppose that $(\rho_n,n\geq 1)$ is a sequence of measures converging
weakly to $\rho$. Assume that $\rho_n$ is $\alpha_n$-s.l.c. for any
$n\geq 1$, then $\rho$ is $\alpha_0=\inf_n\alpha_n$-s.l.c.
\end{lemma}
\begin{lemma}
\label{conv-lemma}
Assume that $d\rho$  is an $\alpha$-s.l.c. measure and
denote by $p_\sigma$ the Gaussian density
$\exp-\frac{1}{2\sigma}|x|^2$.
Then, for any $\delta > 0 $ satisfying 
$$
\frac{1}{\delta}-\frac{1}{\alpha}>\sigma\,,
$$
the measure $\rho\star p_\sigma$   is $\delta$-s.l.c.
\end{lemma}
\nproof
We want to determine the set of $\delta$'s for which the function 
$$
x\to e^{\frac{\delta}{2}|x|^2}\,(\rho\star p_\sigma)(x) 
$$
is log-concave. Let us denote by $\rho_\alpha$ the measure defined by
$$
d\rho_\alpha(y)=\exp\frac{\alpha}{2}|y|^2d\rho(y)\,.
$$
We can write 
\beaa
e^{\frac{\delta}{2}|x|^2}\,(\rho\star p_\sigma)(x)&=&\int
e^{\frac{\delta}{2}|x|^2}
p_\sigma(x-y)e^{-\frac{\alpha}{2}|y|^2}\rho_\alpha(dy)\\
&=&\int
\exp\left[-\half\left|(\frac{1}{\sigma}-\delta)^{1/2}x-\frac{y}{\sigma(\frac{1}{\sigma}-\delta)^{1/2}}\right|^2-\frac{|y|^2}{2}
\left(\frac{\alpha-\delta-\delta\alpha\sigma}{1-\delta\sigma}\right)\right]\rho_\alpha(dy)\,.
\eeaa

Let
$$
p_{\sigma,\alpha}(y)=\exp\left[-\frac{|y|^2}{2}(\frac{\alpha-\delta-\delta\alpha\sigma}{1-\delta\sigma})\right]
$$
which is a Gaussian kernel provided that
$$
\frac{1}{\delta}-\frac{1}{\alpha}>\sigma
$$
and 
$$
x\to e^{\frac{\delta}{2}|x|^2}\,(\rho\star p_\sigma)(x) 
$$
 is a log-concave function from Corollary \ref{cor-2}.
\nqed

\begin{lemma}
If $\rho$ is an $\alpha$-s.l.c. measure on $\R^m$, and if
$F:\R^m\to\R^n$, $m\geq n$, is a linear map,  then $F(\rho)=\rho_F$ is an
$\alpha$-s.l.c. measure.
\end{lemma}
\nproof
Assume first that $\rho$ has a density w.r. to the Lebesgue measure
$l$. We can write $\R^m={\mbox{\rm Im}}(F)+\ker(F)$, then 
$$
\rho_F(f)=\int_{{\mbox{\rm
      Im}}(F)}f(y)\left(\int_{\ker(F)}l(y+y^\bot)dy^\bot
\right)dy\,.
$$
Since, by Corollary \ref{cor-1}, 
$$
y\to\int_{{\mbox{\rm
      Im}}(F)}\exp\frac{\alpha}{2}(|y|^2+|y^\bot|^2)l(y+y^\bot)dy^\bot
$$
is log-concave, 
$$
y\to\exp\frac{\alpha}{2}|y|^2\int_{{\mbox{\rm
      Im}}(F)}l(y+y^\bot)dy^\bot
$$
is also log-concave.
\nqed

\begin{thm}
\label{log-sob1}
Assume that $\rho$ is an $\alpha$-s.l.c. measure on $\R^n$, then it
satisfies the logarithmic Sobolev inequality:
$$
\rho(f^2\log f^2)\leq \frac{2}{\alpha}\rho(|\nabla f|^2)
$$
for any smooth function $f$ with $\rho(f^2)=1$.
\end{thm}
\nproof
Assume first that $d\rho(x)=\rho'(x)dx$, denote by $\mu_\alpha$ the
Gauss measure with covariance $\sqrt{1/\alpha}I_{\R^n}$. Then $\rho\ll
\mu_\alpha$ with
$$
\frac{d\rho}{d\mu_\alpha}(x)=e^{\frac{\alpha}{2}|x|^2}\rho'(x)\,.
$$
By  the hypothesis, this Radon-Nikodym derivative is log-concave, consequently
from a theorem of Caffarelli (cf. \cite{Caf} and \cite{F-U3}), there exists a $1$-Lipschitz map
$T=I_{\R^n}+\nabla \varphi$ such that
$\rho=T(\mu_\alpha)$. Consequently, applying the logarithmic Sobolev
inequality for the Gaussian measure (cf.\cite{Gross})
\beaa
\rho(f^2\log f^2)&=&\mu_\alpha(f^2\circ T \log f^2\circ T)\\
&\leq& \frac{2}{\alpha}\mu_\alpha(|\nabla f\circ
T|^2\|I+\nabla^2\varphi\|^2)\\
&\leq&\frac{2}{\alpha}\mu_\alpha(|\nabla f\circ T|^2)\\
&=&\frac{2}{\alpha}\rho(|\nabla f|^2)\,.
\eeaa
The general case now follows from Lemma \ref{conv-lemma} and a limit
procedure.
\nqed

\section{Infinite dimensional case}
\noindent
The following is basic:
\begin{thm}
\label{inf-dim-thm}
Let $E$ be a separable Fr\'echet space and let $\rho$ be a probability
on $(E,{\calE})$, where $\calE$ denotes the Borel sigma algebra of
$E$. Assume that the finite dimensional projections of $\rho$ are
log-concave. Assume that $f,g,h$ are positive Borel functions
satisfying 
$$
h(su+tv)\geq f(u)^sg(v)^t
$$
for any $u,v\in E$ and $s+t=1$. Then $\rho$ satisfies the
Pr\'ekopa-Leindler property:
$$
\rho(h)\geq \rho(f)^s\rho(g)^t\,.
$$
\end{thm}
\nproof
We can suppose that $\rho$ has convex, compact support $K$ and replace
$E$ by $\R^{\N}$ by injection. Then we can replace $K$ by a product of
compact intervals $J=\prod_{1}^\infty J_n$. If $f$ and $g$ are
continuous, cylindrical functions on $J$, for $x\in J$, define 
$$
k(x)=\sup\left(f(u)^sg(v)^t:\,x=su+tv,\,u,v\in J\right)\,.
$$
The function $k$ is then continuous on $J$ and we have, from the
finite dimensional case, 
$$
\rho(k)\geq \rho(f)^s\rho(g)^t\,.
$$
If $f$ and $g$ are upper semi-continuous on $J$, there exist two
sequences of continuous and cylindrical  functions $(f_n)$ and
$(g_n)$, decreasing to $f$ and $g$ respectively. Hence, $(k_n,\,n\geq
1)$, where $k_n$ is defined as above, converges to $k$ as defined
above and we again have
$$
\rho(k)\geq \rho(f)^s\rho(g)^t\,.
$$
Finally, if $f$ and $g$ are only Borel measurable, there exist two
monotone, increasing sequences $(f_n),\,(g_n)$ whose elements are
upper semi-continuous such that $\lim_n\rho(f_n)=\rho(f)$ and
$\lim_n\rho(g_n)=\rho(g)$. Since we have $h\geq k_n$, it follows that 
\beaa
\rho(h)&\geq&\sup_n\rho(k_n)\geq \sup_n\rho(f_n)^s\rho(g_n)^t\\
&=&\rho(f)^s\rho(g)^t\,.
\eeaa
\nqed

\noindent
The proof of the above theorem contains also the proof of the
following 
\begin{lemma}
\label{contin-lemma}
In order the Pr\'ekopa-Leindler to hold it is necessary and sufficient
that it holds only for the continuous functions $f,g,h$ such that
$h(sx+ty)\geq f(x)^sg(y)^t$ for any $x,y\in E$ and $s+t=1$.
\end{lemma}

\noindent
The proof of the following theorem is quite similar to the proof of
Theorem \ref{inf-dim-thm}:
\begin{thm}
\label{product-thm}
Assume that $E$ and $F$ are two separable Fr\'echet spaces with $\rho$
and $\nu$ satisfying the Pr\'ekopa-Leindler property on $E$ and $F$
respectively. Then the product measure $\rho\otimes \nu$ 
satisfies also  the Pr\'ekopa-Leindler property on $E\times F$.
\end{thm}

\noindent
The following definition is now justified:
\begin{definition}
A Radon  measure on a locally convex space $E$ is called log-concave if it
satisfies the Pr\'ekopa-Leindler property.
\end{definition}

\noindent
The following result is   immediate:
\begin{proposition}
The image of a log-concave measure under any linear, continuous map
is again log-concave.
\end{proposition}

\noindent
From Lemma \ref{contin-lemma} we get at once
\begin{corollary}
Let $(\rho_n,n\geq 1)$ be a sequence of log-concave measures
converging weakly to a measure $\rho$, then $\rho$ is also
log-concave.
\end{corollary}

\begin{corollary}
Let $\rho$ be a bounded measure on $\R^{\N}$, let us denote by
$(\pi_n,n\geq 1)$ the canonical finite dimensional projections. The measures
$(\pi_n(\rho),n\geq 1)$ are  log-concave if and only if  $\rho$ is
log-concave.
\end{corollary}
\nproof
We can write $\rho$ as the weak limit of the sequence of measures
$(\pi_n(\rho)\otimes \delta^n,n\geq 1)$ where $\delta^n$ is the image under $I-\pi_n$ of the
Dirac $\delta_0$ measure on $\R^\N$.
\nqed

\noindent
Assume now that $E$ is a separable Fr\'echet space and assume that $H$
is a separable  Hilbert space continuously and densely injected into $E$. We
identify $H$ with its continuous dual, hence $E^\star\subset H\subset
E$. Choose a sequence $(\tilde{e}_i,i\geq 1)$ from $E^\star$ such that
its image under the injection $E^\star\hookrightarrow H$, denoted as
$(e_i,i\geq 1)$ is a complete, orthonormal base of $H$. Define $\pi_n$
on $E$ as 
$$
\pi_n(x)=\sum_{i\leq n}\langle x,\tilde{e}_i\rangle e_i\,.
$$
The typical examples  for this situation is the case of the Wiener
space for $E$ and the Cameron-Martin space for $H$  or $E=\R^\N,
\,H=l^2$. 

\begin{definition}
\label{slc-inf}
A  Radon measure $\rho$  on $E$ is called $\alpha$-s.l.c. 
\begin{enumerate}
\item if $\lim_n\pi_n(x)=x$ $\rho$-almost everywhere,
\item if $\pi_n(\rho)$ is $\alpha$-s.l.c. on the Euclidean space
  spanned by $\{e_1,\ldots,e_n\}$, for any $n\geq 1$.
\end{enumerate}
\end{definition}

\noindent
We have the following result which is the immediate consequence of the
finite dimensional case (cf. Theorem \ref{log-sob1}):
\begin{thm}
\label{log-sob2}
If $\rho$ is an $\alpha$-s.l.c. measure on $E$, then it satisfies the
logarithmic Sobolev inequality:
$$
\rho(f^2\log f^2)\leq \frac{2}{\alpha}\rho(|\nabla f|_H^2)
$$
for any smooth, cylindrical function $f$ with $\rho(f^2)=1$.
\end{thm}

\section{The case of abstract Wiener space}
\noindent
While working in this frame one encounters often the difficulty of
defining a proper Jacobian due to the lack of regularity of the
corresponding transformation. This happens especially in the case of
the measure transportation theory.  Consequently it is reasonable to push the ways to extend
as much as possible the notion of Jacobian of a transformation with
unsufficient regularity. 
\subsection{Sub-jacobians for monotone transformations}
Let $(W,H,\mu)$ be an abstract Wiener space, we say that a map
$U=I_W+u$, where $u:W\to H$ is a measurable map is monotone or a
monotone shift,  if $h\to
(h+u(w+h),h)_H\geq 0$ $\mu$-almost surely.
\begin{lemma}
\label{lemma-1} Assume that $U=I_W+u$ is a monotone shift with
$u\in \DD_{p,1}(H),\,p>1$. Then
$$
E[f\circ U\,\La(U)]\leq E[f]\,,
$$
for any positive $f\in C_b(W)$, where
$$
\La(U)=\dett(I+\nabla u)\exp\left(-\delta
u-\frac{1}{2}|u|_H^2\right)\,
$$
and $\dett(I_H+\nabla u)$ denotes the modified Carleman-Fredholm determinant.
\end{lemma}
\proof The necessary background about the subject and the  proof
follows from \cite{BOOK}, Theorem 6.3.1. 
\nqed

\remark
If $A$ is a nuclear operator on a Hilbert space, then
$\dett(I_H+A)$ is defined as 
$$
\dett(I_H+A)=\det(I_H+A)\exp-\trace(A)\,,
$$
and this function has an analytic extension to the space of
Hilbert-Schmidt operators, consequently, the log-concavity of the
ordinary determinant implies the log-concavity of the map
$A\to \dett(I_H+A)$.

\begin{lemma}
\label{lemma-2} Assume that $U=I_W+u$ is a monotone shift and
$u\in\DD_{p,0}(H)=L^p(\mu,H)$ for some $p>1$. Then there exists
some $\la(u)\geq 0$ a.s., $E[\la(u)]\leq 1$ and
$$
E[f\circ U\,\la(u)]\leq E[f]
$$
for any positive and measurable $f$. In particular, $U\mu$ is
absolutely continuous w.r.to $\mu$ on the set on which $\la(u)>0$.
\end{lemma}
\proof Let $(P_t,t\geq 0)$ be the Ornstein-Uhlenbeck semigroup,
let $U_n=I_W+u_n$, with $u_n=P_{1/n}u$. Then apply Lemma
\ref{lemma-1} to $U_n$. Define
$$
\la(U)=\lim\inf_n\La(U_n)\,.
$$
$\la(U)$ does exist and, from the Fatou lemma,  it is integrable
and satisfies the claim.
\nqed
\begin{corollary}
\label{cor-11} 
The map $U\to \la(U)$ is almost surely log-concave
on the set of monotone shifts.
\end{corollary}
\proof
If $U_i=I_W+u_i$, with $u_i\in L^p(\mu,H)$,  $i=1,2$ are
two monotone shifts, define $u_i^n=P_{1/n}u_i$, $U_i^n=I_W+u_i^n$,
$i=1,2$. Then, by the log-concavity of the Carleman-Fredholm
determinant,  for $s,\,t\in [0,1]$ with $s+t=1$, we have 
$$
\La(sU_1^n+tU_2^n)\geq \La(U_1^n)^s\,\La(U_2^n)^t
$$
a.s. If we take the $\lim\inf$ of both sides as $n\to \infty$, the
inequality is preserved.
\nqed
\begin{proposition}
\label{prop-11} Assume that $U_n=I_W+u_n$, $n\geq 1$, is a sequence
of shifts, where $u_n\in\DD_{p,1}(H)$, for some $p>1$. Assume that
$U_n\mu\ll\mu$ for all $n\geq 1$ and denote
$$
L_n=\frac{dU_n\mu}{d\mu}\,.
$$
Assume further that $L_n\circ U_n\,\La(U_n)=1$ a.s. for all $n\geq
1$, also that $L_n\to L$ in $L^1(\mu)$ and finally that $U_n\to
U=I_W+u$ in $L^0(\mu,W)$. Then $L_n\circ U_n\to L\circ U$,
$\La(U_n)\to \la(U)$ in $L^0(\mu)$ and we have
$$
L\circ U\,\la(U)=1
$$
a.s.
Besides, under the additional hypothesis:
$$
\sup_nE[\log^+ L_n]<\infty\,,
$$
the sequence $(\La(U_n),n\geq 1)$ is uniformly integrable, hence it
converges to $\la(U)$ also in $L^1(\mu)$ and we have 
$$
E[f\circ U\,\la(U)]=E[f]\,,
$$
for any $f\in C_b(W)$.
\end{proposition}
\proof Since $L_n\to L$ in $L^1(\mu)$ and $U_n\to U$ in
probability, it follows from the Lusin theorem that $L_n\circ
U_n\to L\circ U$ in probability, hence $(\La(U_n),\,n\geq 1)$
converges in probability and hence
$$
\lim\La(U_n)=\lim\inf_n\La(U_n)=\la(U)\,.
$$
The rest is obvious, since the last hypothesis implies precisely the
uniform integrability of the sequence $(\La(U_n),\,n\geq 1)$.
\nqed
\begin{remarkk}
Proposition \ref{prop-11} is astonishing in the sense that we do not make
any assumption about the convergence of the $H$-valued parts of the
shifts at all. Assume now  that $u_n\to u$ in
$L^0(\mu,H)$, then, although $\delta u$ either $\nabla u$ do not
exist, the sequence
$$
(\dett(I+\nabla u_n)e^{-\delta u_n},\,n\geq 1)
$$
converges in probability to a non-trivial limit that we denote by
$J(U)$, hence $\la(U)$ can be represented as
$$
\la(U)=J(U)\exp-\frac{1}{2}|u|^2\,.
$$
\end{remarkk}

\subsection{The transport case and applications}
Assume that $L\in L\log L(\mu)$, $E[L]=1$. Let $T=I_W+\nabla \phi$,
with  $\phi$ in $\DD_{2,1}$,
be the transport map which maps $d\mu$ to $Ld\mu$ whose properties are
proved in \cite{F-U2}.
Define $L_n=E[P_{1/n}L|V_n]$, where $V_n$ is the sigma algebra
generated by $\{\delta e_1,\ldots,\delta e_n\}$,
$(e_i,i\geq)\subset W^\star$ being an orthonormal basis of $H$.
Let $T_n=I_W+\nabla\phi_n$, $\phi_n\in \DD_{2,1}$ be the transport
map which maps $d\mu$ to $L_nd\mu$. Recall that $\phi_n$ is
$1$-convex, since it is $V_n$-measurable, $\nabla^2\phi$ and
$\calL\phi$ are lower bounded distributions, hence they are measures. We have
$$
E[f\circ T_n]=E[f\,L_n]
$$
for any $f\in C_b(W)$. Besides, since $L_n>0$ a.s., we have
$$
L_n\circ T_n\,\La(T_n)=1
$$
a.s., where
$$
\La(T_n)=\dett(I+\nabla_a^2\phi_n)\exp\left(-\calL_a\phi_n-\frac{1}{2}|\nabla\phi_n|^2\right)\,,
$$
where $\nabla_a^2\phi_n$ and $\calL_a\phi_n$ denote respectively the
absolutely continuous parts of the measures $\nabla^2\phi_n$ and
$\calL\phi_n$ and $\dett$ denotes, as usual, the modified
Carleman-Fredholm determinant (c.f.\cite{BOOK} for further information).
It follows from \cite{F-U2} that $\phi_n\to \phi$ in $\DD_{2,1}$,
$(L_n,n\geq 1)$ being uniformly integrable, $L_n\circ T_n\to
L\circ T$ in probability, hence $\La(T_n)\to \la(T)$ also in
probability, where $\la(T)$ can be represented as
$$
\la(T)=J(T)\exp-\frac{1}{2}|\nabla \phi|^2\,.
$$
In \cite{F-U2}, we have shown that the sequence $(\calL_a\phi_n,n\geq
1)$ is a submartingale with respect to the increasing sequence of
sigma algebras $(V_n,n\geq 1)$ and  the  inequality $\calL_a\phi_n\leq -\log\La_n$
implies 
$$
(\calL_a\phi_n)^+\leq (-\log\La(T_n))^+\,.
$$
Consequently 
\beaa
E[(\calL_a\phi_n)^+]&\leq&E[(-\log\La(T_n))^+]=E[(\log L_n\circ T_n)^+]\\
&=&E[L_n\log^+L_n]\leq 2e^{-1}+E[L_n\log L_n]
\eeaa
and Jensen inequality implies that 
$$
\sup_nE[(\calL_a\phi_n)^+]<\infty\,,
$$
which is a sufficient condition for the almost everywhere convergence
of the submartingale $(\calL_a\phi_n,n\geq 1)$ whose limit we denote
as $\calL(\phi)\in L^1(\mu)$. Note that, as a consequence of this
observation, combined with the convergence of $(\La(T_n),n\geq 1)$, we
deduce also the convergence of $(\dett(I_H+\nabla^2_a\phi_n),n\geq 1)$
in probability.
Moreover, we have also
$$
E[f\circ T\,\la(T)]\leq E[f]\,,
$$
for any measurable, positive $f$ and moreover, this inequality becomes
an equality if $\log L$ is integrable.

\begin{proposition}
\label{prop-2} Let $L_1,\,L_2$ be as above, denote by $T_1,\,T_2$
the corresponding transport maps and define $M=aT_1+bT_2$, where
$a+b=1, a\geq 0$. Define $M_n=a T_1^n+b T_2^n$ as in Corollary
\ref{cor-11} and define finally $\la(M)=\lim\inf_n\La(M_n)$. We
then have $M\mu$ is absolutely continuous w.r.to $\mu$ and
$$
\la(M)\geq \la(T_1)^a\,\la(T_2)^b
$$
a.s.
\end{proposition}
\proof We have
\beaa
\la(M)&=&\lim\inf\La(M_n)\\
&\geq&\lim\inf\La(T_1^n)^a\La(T_2^n)^b\\
&\geq&\lim\inf \La(T_1^n)^a\,\lim\inf\La(T_2^n)^b\\
&=&\la(T_1)^a\la(T_2)^b
\eeaa
a.s. Hence $\la(M)>0$ a.s. This result combined with the following
consequence of  the Fatou lemma
$$
\int f\circ M\la(M)d\mu\leq \int fd\mu\,,
$$
for any $0\leq f\in C_b(W)$, implies the absolute continuity $M\mu\ll\mu$.
\nqed
\begin{thm}
\label{thm-1} Assume that $a,b$ and $c$ are measurable, positive
functions on $W$ such that, for given $s,t\in [0,1]$ with $s+t=1$
and for any $h,k\in H$, we have
$$
a(w+sh+tk)\exp\left[-\frac{1}{2}|sh+tk|_H^2\right]\geq
\left(b(w+h)\exp-\frac{1}{2}|h|_H^2\right)^s\,\left(c(w+k)\exp-\frac{1}{2}|k|_H^2\right)^t
$$
almost surely. Let also $q$ be any $H$-logconcave density and
denote by $\nu$ the measure $d\nu=qd\mu$. Then we have
$$
\int a\,d\nu\geq \left(\int b\,d\nu\right)^s\left(\int
c\,d\nu\right)^t\,.
$$
\end{thm}
\nproof First we shall prove the case $q=1$, then the general case
can be reduced to this particular case by replacing $a,b$ and $c$
by $a\,q,\,b\,q$ and by $c\,q$ respectively. Moreover, by
replacing $a,b,c$ by $a\wedge n,\,b\wedge n$ and $c\wedge n$ we
may suppose that they are bounded. Finally, by multiplying them
with adequate constants, we can also suppose that their integrals
w.r.to $\mu$ are all equal to unity. Let $T_1=I+\nabla \phi_1$ and
$T_2=I+\nabla \phi_2$ be the transport maps such that
$T_1\mu=b\cdot\mu,\,T_2\mu=c\cdot\mu$, where $l\cdot\mu$ denotes
the measure with density $l$. It follows from above explanations,
$\la(T_1)$ and $\la(T_2)$ are well-defined and
$$
b\circ T_1\,\la(T_1)=c\circ T_2\,\la(T_2)=1
$$
a.s. Let $M=sT_1+tT_2$, then $M\mu\ll \mu$ and as explained above
$$
\la(M)\geq \la(T_1)^s\,\la(T_2)^t
$$
a.s. Hence
\beaa
1&=&(b\circ T_1)^s(c\circ T_2)^t\la(T_1)^s\la(T_2)^t\\
&\leq&a\circ(sT_1+sT_2)\la(T_1)^s\la(T_2)^t\\
&\leq& a\circ M\,\la(M)\,.
\eeaa
Therefore
$$
1\leq \int a\circ M\,\,\la(M)d\mu\leq \int a\,\, d\mu
$$
and this accomplishes the proof.
\nqed

\remark
Here is another proof of the theorem: let $(\pi_n,n\geq 1)$ be a
sequence of orthogonal projections, constructed from the elements of
$W^\star$,  of the Cameron-Martin space $H$
increasing to identity such that $\lim \pi_nw=w$ $\mu$-a.s. Let
$w_n=\pi_nw$ and $w_n^\bot=w-w_n$. For a measurable function $f$ on
$W$, denote the partial map $w_n\to f(w_n+w_n^\bot)\exp-\frac{1}{2}|w_n|^2$ by
$f_{w_n^\bot}$. The hypothesis above is equivalent to (cf.\cite{F-U1})
$$
a_{w_n^\bot}(sx+ty)\geq b_{w_n^\bot}(x)^s\, c_{w_n^\bot}(y)^t\,.
$$
Since $\mu_n=\pi_n\mu$ is log-concave, we have 
$$
E[a|\pi_n^\bot]\geq E[b|\pi_n^\bot]^s E[c|\pi_n^\bot]^t
$$
and the (second) proof follows from the (reverse) martingale
convergence theorem.

\noindent
The following concept has been studied already in \cite{F-U1}:

\begin{definition}
A measurable, $\R_+$-valued function on $W$ is called $1$-log
concave if, for any $s+t=1,\, s\geq 0$, \beaa
\lefteqn{f(w+sh+th')\exp\left(-\frac{1}{2}|sh+th'|_H^2\right)}\\
&\geq&
\left(f(w+h)\exp-\frac{1}{2}|h|_H^2\right)^s\left(f(w+h')\exp-\frac{1}{2}|h'|_H^2\right)^t
\eeaa
 almost surely for any $h,h'\in H$.
\end{definition}

\begin{corollary}
\label{cor-12}
 Assume that $(W_1,H_1,\mu_1)$ and $(W_2,H_2,\mu_2)$
be two abstract Wiener spaces. Assume that $f:W_1\times W_2\to
\R_+$ is an $1$-log concave on the abstract Wiener space
$(W_1\times W_2,H_1\times H_2,\mu_1\times \mu_2)$. Then
$$
\hat{f}(x)=\int_{W_2}f(x,y)\,\mu_2(dy)
$$
is $1$-log concave on $(W_1,H_1,\mu_1)$.
\end{corollary}
\nproof Let $h,h'\in H_1$,  $s+t=1$. For $(x,y)\in W_1\times W_2$,
define
$$
a_{x+sh+th'}(y)=f((x,y)+s(h,0)+t(h',0))\exp\left(-\half|sh+th'|_{H_1}^2\right)
    \,.
$$
For any $k,k'\in H_2$, we have
 \beaa
\lefteqn{a_{x+sh+th'}(y+sk+tk')\exp\left(-\half|sk+tk'|_{H_2}^2\right)}\\
& \geq& a_{x+h}(y+k)^s
\exp\left(-\frac{s}{2}|k|_{H_2}^2\right)\,a_{x+h'}(y+k')^t\left(\exp-\frac{t}{2}|k'|_{H_2}^2\right)
\,.
\eeaa
 Applying Theorem \ref{thm-1} to $ a_{x+sh+th'},\,a_{x+h}$
and $a_{x+h'}$, we get
\begin{gather*}
\hat{f}(x+sh+th')\exp\left(-\half|sh+th'|_{H_1}^2\right)=\int_{W_2}a_{x+sh+th'}(y)\mu_2(dy)\\
\geq\left(\int_{W_2}a_{x+h}(y)\mu_2(dy)\right)^s\left(\int_{W_2}a_{x+h'}(y)\mu_2(dy)\right)^t\\
=\hat{f}(x+h)^s\hat{f}(x+h')^t\exp-\frac{s}{2}|h|_{H_1}^2\exp-\frac{t}{2}|h'|_{H_1}^2\,.
\end{gather*}
\nqed

\begin{corollary}
\label{cor-14} Assume $B,C$ are two measurable subsets of $W$,
then, for any measure $\nu$ as in Theorem \ref{thm-1}, we have
$$
\nu(sB+tC)\geq \nu(B)^s\nu(C)^t\,,
$$
where $s+t=1,\, s,\,t\geq 0$.
\end{corollary}
\nproof Let $a,\,b$ and $c$ be the indicator functions of the sets
$sB+tC$, $B$ and $C$ respectively. \nqed

\noindent
 Before further infinite dimensional  considerations let
us  give a lemma which is a version of Theorem 2.1 of \cite{Bor-1}
\begin{lemma}
\label{std-lemma} Assume that $f,f_0,f_1$ are positive, bounded,
measurable functions on $\R^n$ such that for any Borel sets
$C_0,C_1$, some  $\al,\beta\in [0,1]$ with $\al+\beta=1$, we have
\begin{equation}
\label{inek}
\int_{C}f(x)dx\geq
\left(\int_{C_0}f_0(x)dx\right)^\al\left(\int_{C_1}f_1(x)dx\right)^\beta\,,
\end{equation}
where $C=\al C_0+\beta C_1$. Then, for any $x_0,x_1\in \R^n$, we
have
$$
f(z+\al x_0+\beta x_1)\geq f_0(z+x_0)^\al\,f_1(z+x_1)^\beta
$$
$dz$-almost surely. In particular, if the above identity holds
whenever we replace the  Lebesgue integral with the Gaussian
integral, then $f, f_0$ and $f_1$ satisfy the following identity
\begin{gather*}
f(z+\al x_0+\beta x_1)\exp\half|\al x_0+\beta x_1|^2\\
\geq\left(f_0(z+x_0)\exp-\frac{1}{2}|x_0|^2\right)^\al\,\left(f_1(z+x_1)
\exp-\frac{1}{2}|x_1|^2\right)^\beta
\end{gather*}
$dz$-almost surely.
\end{lemma}
\nproof We may suppose that the functions are of bounded support.
Let $C_0=z+x_0+\eps I_n$, $C_1=z+x_1+\eps I_n$, where
$I_n=[-1/2,1/2]^n$. Then the l.h.s. of  the  inequality (\ref{inek}) can be
written as
$$
\frac{1}{\eps^n}\int 1_{\eps I_n}(x-(z+\al x_0+\beta x_1))f(x)dx
$$
which converges in measure to $f(z+\al x_0+\beta x_1)$ as $\eps\to
0$.  For the terms at the r.h.s. we have similar convergence
results in measure. For the Gaussian case, it suffice to replace
the functions $f,f_0,f_1$ with $f\,q,f_0\,q$ and $f_1\,q$
respectively, where $q$ denotes the Gaussian density.
 \nqed

Let $(e_n,n\geq 1)\subset W^*$ be a CONB basis of $H$ and denote
by $V_n$ the sigma algebra generated by $\{\delta
e_1,\ldots,\delta e_n\}$ and completed with $\mu$-negligeable
sets. We have
\begin{thm}
\label{thm-12} Assume that $a,b$ and $c$ are measurable, positive
functions on $W$. For given $s,t\in [0,1]$ with $s+t=1$ and for
any $h,k\in H$, we have
\begin{enumerate}
\item For any $h,k\in H$, the following inequality holds
$\mu$-a.s.
\begin{gather}
\label{prekopa}
 a(w+sh+tk)\exp-\half|sh+tk|_H^2\geq
 \left(b(w+h)\exp-\half|h|_H^2\right)^s\nonumber\\
 \left(c(w+k)
 \exp-\half|k|_H^2\right)^t
\end{gather}
 if and only if
\begin{gather*}
E[a|V_n](w+s\pi_nh+t\pi_nk)\exp-\half|s\pi_nh+t\pi_nk|_H^2\\
\geq\left(E[b|V_n](w+\pi_nh)\exp-\half|\pi_nh|_H^2\right)^s\\
\cdot\left(E[c|V_n](w+\pi_nk)\exp-\half|\pi_nk|_H^2\right)^t
\end{gather*}
$\mu$-a.s., where $\pi_n$ denotes the orthogonal projection of $H$
onto the space spanned by $\{e_1,\ldots,e_n\}$. \item Similarly,
the relation (\ref{prekopa}) is equivalent to
\beaa
 P_\tau a(w+sh+tk)\exp-\half|sh+tk|_H^2&\geq&
 \left(P_\tau b(w+h)\exp-\half|h|_H^2\right)^s\\
&& \left(P_\tau c(w+k)\exp-\half|k|_H^2\right)^t
\eeaa
 for any $\tau\geq 0$, where $P_\tau$ denotes the Ornstein-Uhlenbeck
semigroup.
\end{enumerate}
\end{thm}
\nproof We can suppose that $a,b,c$ are bounded. Denote by
$a_n,b_n,c_n$ the conditional expectations of $a,b,c$ respectively
w.r.to $V_n$. Let now $B$ and $C$ be $V_n$-measurable sets, hence
$sB+tC$ is also $V_n$-measurable (recall that $V_n$ is
completed!). It then follows from Theorem \ref{thm-1}
\beaa
\int_{sB+tC}a_n\,d\mu&=&\int_{sB+tC}a\,d\mu\\
&\geq&\left(\int_B b\,d\mu\right)^s\left(\int_C c\,d\mu\right)^t\\
 &=&\left(\int_B b_n\,d\mu\right)^s\left(\int_C
 c_n\,d\mu\right)^t\,.
\eeaa Since this true for any $V_n$-measurable set, it follows
from  Lemma \ref{std-lemma} that $a_n,b_n$ and $c_n$ satisfy the
inequality claimed in the first part of the theorem (with the
Gaussian measure). To prove the second part, we can replace
$a,b,c$ by $a_n,b_n,c_n$ since $P_\tau$ commutes with the
conditional expectation w.r.to $V_n$. Hence the problem is reduced
to the finite dimensional case. Let us denote again by the same
notation the Ornstein-Uhlenbeck semigroup on $\R^n$. Then we can
write
$$
\int_{sB+tC}P_\tau
a_n(x)d\mu(x)=\int_{sB+tC}\int_{\R^n}a_n(y)q_\tau(x,y)dy
d\mu(x)\,,
$$
where
$$
q_\tau(x,y)=(2\pi(1-e^{-2\tau}))^{-n/2}\exp-\frac{|y-e^{-\tau}x|^2}{2(1-e^{-2\tau})}\,,
$$
which is a log-concave function in two variables. We can also
write
$$
(sB+tC)\times \R^n=s(B\times \R^n)+t(C\times \R^n)\,.
$$
It then follows from the Pr\'ekopa-Leindler inequality in
$\R^n\times \R^n$ that
$$
\int_{sB+tC}P_\tau a_n d\mu \geq \left(\int_{B}P_\tau b_n
d\mu\right)^s \left(\int_{C}P_\tau c_n d\mu\right)^t\,.
$$
It follows then from Lemma \ref{std-lemma} that we have \beaa
P_\tau a_n(w+sh+tk)\exp-\half|s\pi_nh+t\pi_nk|_H^2&\geq&
\left(P_\tau
b_n(w+h)\exp-\half|\pi_nh|_H^2\right)^s\\
&& \left(P_\tau c_n(w+k)\exp-\half|\pi_nk|_H^2\right)^t
\eeaa
 almost surely
and we can pass to the limit as $n\to\infty$ due to the martingale
convergence theorem. \nqed



{\footnotesize}

\end{document}